\documentclass[12pt]{article}
\usepackage{amsmath,amsfonts}
\usepackage{amssymb}
\usepackage{amsthm}
\usepackage{times,txfonts}

\usepackage{accents, kbordermatrix}

\usepackage{pict2e}

\newtheorem{thm}{Theorem}[section]
\newtheorem{prop}[thm]{Proposition}

\newtheorem{lemma}[thm]{Lemma}

\newtheorem{remark}[thm]{\it Remark}

\numberwithin{equation}{section}

\def\pf{\noindent{\it Proof.} \ }

\def\qed{\hfill $\square$}

\def\dim{{\rm dim}~}

\def\hsymb#1{\mbox{\strut\rlap{\smash{\huge$#1$}}\quad}}

\title{Hypergeometric solution of a certain
polynomial Hamiltonian system 
of isomonodromy type}

\author{Teruhisa Tsuda \\
Faculty of Mathematics, Kyushu University,  
\\
Fukuoka 819-0395, Japan}

\date{April 30, 2010 (Revised: September 13, 2010)}

\begin{document}
\maketitle

\begin{abstract} 
In \cite{tsu10}
a unified description as polynomial Hamiltonian systems 
was
established for 
a broad class of the Schlesinger systems
including the sixth Painlev\'e equation
and Garnier systems.
The main purpose of this paper is to present 
particular solutions of
this Hamiltonian system
in terms of a certain generalization of Gau\ss's hypergeometric function.
Key ingredients of the argument  
are the linear Pfaffian system derived from 
an integral representation of the hypergeometric function
(with the aid of twisted de Rham theory)
and Lax formalism of the Hamiltonian system.
\end{abstract}

\renewcommand{\thefootnote}{\fnsymbol{footnote}}
\footnotetext{{\it 2000 Mathematics Subject Classification} 
33C70, 
34M55, 
37K10. 
} 
\footnotetext{{\it Keywords:} 
isomonodromic deformation, 
hypergeometric function, 
Painlev\'e equation,
Schlesinger system.} 
\footnotetext{{\it E-mail:} 
 tudateru@math.kyushu-u.ac.jp}

\section{Introduction}
Fix  integers $L \geq 2$ and $N \geq 1$.
We consider
 the following completely integrable
 Hamiltonian system 
 of partial differential equations:
\begin{subequations} \label{subeq:HLN}
\begin{equation} \label{eq:HLN_caneq}
\frac{\partial q_n^{(i)}}{\partial x_j}=\frac{\partial H_j}{\partial p_n^{(i)}},
\quad 
\frac{\partial p_n^{(i)}}{\partial x_j}=-\frac{\partial H_j}{\partial q_n^{(i)}}
\quad 
\left(
\begin{array}{l}
  i,j =1, \ldots, N\\
 n=1, \ldots,   L-1
\end{array}
\right)
\end{equation}
with variables $x=(x_1,\ldots,x_N)$
and unknowns $q_n^{(i)}$ and $p_n^{(i)}$.
Here the Hamiltonians function $H_i$ is given by 
\begin{equation} \label{eq:HLN_ham}
x_i H_i=
\sum_{n=0}^{L-1} e_n q_n^{(i)} p_n^{(i)} 
+\sum_{j=0}^N \sum_{0 \leq m <n \leq L-1}
 q_m^{(i)}   p_m^{(j)} q_n^{(j)} p_n^{(i)} 
+\sum^N_{
\begin{subarray}{c} 
j=0  \\
j \neq i
\end{subarray}
}
\frac{x_j}{x_i-x_j} 
 \sum_{m,n=0}^{L-1}
 q_m^{(i)}   p_m^{(j)} q_n^{(j)} p_n^{(i)}
\end{equation}
\end{subequations}
and
\[
x_0=q_n^{(0)}
= q_0^{(i)}=1, \quad
p_n^{(0)}
=
\kappa_{n} - \sum_{i=1}^{N} q_n^{(i)} p_n^{(i)},
\quad
p_0^{(i)}
=\theta_i-\sum_{n=1}^{L-1} q_n^{(i)} p_n^{(i)},
\]
thereby, $H_i$ forms a polynomial in the canonical 
variables.
This system
was introduced 
by the author 
in \cite{tsu10} 
 via a similarity reduction of 
 the UC hierarchy \cite{tsu04}
(an extension of the KP hierarchy);
it is equivalent to a class of the Schlesinger systems
\cite{sch12}
describing isomonodromic deformations
of an $L \times L$
Fuchsian system with $N+3$ poles on the Riemann sphere.
The {\it spectral type} of this Fuchsian system is given by
the $(N+3)$-tuple 
\[
(1,1,\ldots,1),(1,1,\ldots,1),
 \underbrace{(L-1,1), \ldots, (L-1,1)}_{N+1} 
\]
of partitions of $L$,
which
indicates how the characteristic exponents overlap 
at each of the $N+3$ regular singularities.
System (\ref{subeq:HLN})
contains complex constants
\[
(e,\kappa,\theta)=(e_0,e_1,\ldots,e_{L-1},\kappa_0,\kappa_1,\ldots,\kappa_{L-1},\theta_0,\theta_1,\ldots,\theta_N)
\]
satisfying
the linear constraints 
\[
\sum_{n=0}^{L-1}e_n=\frac{L-1}{2}
\quad
\text{and} \quad
\sum_{n=0}^{L-1}\kappa_n
=\sum_{i=0}^N \theta_i,
\]
so the number of constant 
parameters is essentially $2L+N-1$.
Note that these parameters
correspond to the characteristic exponents of the associated
Fuchsian system.
Henceforth, we denote 
by
${\cal H}_{L,N}$ 
the polynomial Hamiltonian system
(\ref{subeq:HLN}).
For example, 
the case where $L=2$ and any $N \geq 1$ 
(${\cal H}_{2,N}$) 
coincides with 
the Garnier system in $N$ variables \cite{gar12, ko84}
and, thus, the first nontrivial case
(${\cal H}_{2,1}$)
with the sixth Painlev\'e equation
$P_{\rm VI}$
\cite{mal22, oka87}.
For details refer to \cite{tsu10}.

In this paper  we present
a family of particular solutions of 
the polynomial Hamiltonian system
${\cal H}_{L,N}$,
which is
parameterized by 
a point in the projective space ${\mathbb P}^{N(L-1)}$. 
These solutions are governed by a linear Pfaffian system of rank 
$N(L-1)+1$ and, furthermore, expressed in terms of a 
generalization of Gau\ss's hypergeometric function.

We begin by introducing 
the hypergeometric function
crucial to this work.
Fix the notation of multi-index;
let $I=\{m=(m_1,\ldots, m_N) \ | \ m_i \in {\mathbb Z}_{\geq 0} \}$, and 
write
$x^m={x_1}^{m_1} \cdots {x_N}^{m_N}$
and $|m|=m_1+ \cdots +m_N$ for $m \in I$.
We define a function 
$F_{L,N}=F_{L,N}(\alpha, \beta, \gamma; x)$
in $N$ variables $x=(x_1,\ldots,x_N)$
by means of the power series
\begin{equation} \label{eq:FLN}
F_{L,N}(\alpha, \beta, \gamma; x)
\stackrel{\rm def}{=} 
\sum_{m \in I}
\frac{(\alpha_1)_{|m|} \cdots (\alpha_{L-1})_{|m|} (\beta_1)_{m_1} \cdots  (\beta_N)_{m_N} }{ (\gamma_1)_{|m|} \cdots (\gamma_{L-1})_{|m|} (1)_{m_1} \cdots (1)_{m_N}} 
x^m
\end{equation}
convergent in the polydisc $\{|x_1|<1, \ldots,|x_N| <1\} \subset {\mathbb C}^N$.
Here 
\[
(\alpha,\beta,\gamma)=(\alpha_1,\ldots,\alpha_{L-1}, \beta_1,\ldots,\beta_N, \gamma_1, \ldots,\gamma_{L-1})
\]
are complex constants such that 
$\gamma_n \notin {\mathbb Z}_{< 0}$,
and 
$(a)_n=\Gamma(a+n)/\Gamma(a)$.
The series $F_{L,N}$ satisfies the system 
of linear differential equations
\begin{equation} \label{eq:hgeq}
\left\{
x_i  \left( \beta_i+\delta_i\right)
\prod_{k=1}^{L-1} \left( \alpha_k+ 
{\cal D}
\right) 
- \delta_i
 \prod_{k=1}^{L-1} \left( \gamma_k-1+  {\cal D} \right) 
\right\}y=0,
\quad i=1,\ldots,N,
\end{equation}
where 
\begin{equation} \label{eq:eular}
\delta_i=x_i \frac{\partial}{\partial x_i}
\quad \text{and} \quad
{\cal D}=\sum_{i=1}^N \delta_i.
\end{equation}
As is shown in the next section,
system (\ref{eq:hgeq}) is equivalent to a linear
Pfaffian system of rank $N(L-1)+1$.
The holomorphic function $F_{L,N}(\alpha, \beta, \gamma; x)$ 
at
$0 \in {\mathbb C}^N$
can be analytically continued 
along any path on 
$X=\{ 
x=(x_1,\ldots,x_N) \in {\mathbb C}^N \ | \  
x_i \neq x_j  \ (i \neq j), \
x_i \neq 0,1
\}$,
which is the complement of the singular locus of this Pfaffian system.
Note that,
if $(L,N)=(2,1)$, $(L,1)$
and
$(2,N)$, then
the hypergeometric series $F_{L,N}$ reduces  
to Gau\ss's ${}_2F_1$, 
Thomae's ${}_L F_{L-1}$ \cite{tho70} 
and 
Appell--Lauricella's $F_D$ \cite{ak26, lau93},
respectively.

We turn now to the particular solutions of the Hamiltonian system 
${\cal H}_{L,N}$.
There is an equivalent formulation of ${\cal H}_{L,N}$
(Lax formalism)
as a compatibility condition of 
an auxiliary linear problem 
consisting of the foregoing Fuchsian system 
and its deformation equations;
see \cite{tsu10}.
Under a certain condition of parameters
we find particular solutions of  ${\cal H}_{L,N}$
such that the associated Fuchsian system 
becomes reducible;
furthermore, these solutions 
are governed by
the same linear Pfaffian system of rank $N(L-1)+1$ 
as $F_{L,N}$.
This fact leads us to the

\begin{thm} \label{thm:main}
When $\kappa_0-\sum_{i=1}^N \theta_i=0$,
 the Hamiltonian system ${\cal H}_{L,N}$ 
possesses an $N(L-1)$-parameter family of particular solutions,
each of which is
expressed in terms of a hypergeometric function, i.e., 
an arbitrary solution of  
{\rm(\ref{eq:hgeq})}.
\end{thm}
\noindent
(See Theorem~\ref{thm:hgsol}.)

\begin{remark}\rm
The hypergeometric solutions of 
$P_{\rm VI}$ $(={\cal H}_{2,1})$ and the Garnier system 
$(={\cal H}_{2,N})$
were
first
given by Fuchs \cite{fuc07} 
and Okamoto--Kimura \cite{ok86},
respectively.
See also \cite{gar12}.
They linearized
the Riccati-type equations 
for particular solutions
by introducing new dependent variables,
and then identified the resulting linear ones as the hypergeometric 
differential equations.
Recently the case of ${\cal H}_{L,1}$ 
was studied independently by Suzuki \cite{suz10}; 
he obtained a power series solution
through Frobenius' method
after a direct linearization.
The present result covers all previous ones;
however, it is based on a method quite different from theirs.
We emphasize that 
key ingredients of
the argument are 
a systematic derivation of 
the linear Pfaffian system 
from an integral representation of the hypergeometric function
(with the aid of twisted de Rham theory)
and investigation into its Lax formalism 
rather than
the Hamiltonian system itself.
\end{remark}

In Sect.~\ref{sect:hg} 
we present the integral representation of the hypergeometric function
$F_{L,N}$ 
(Proposition~\ref{prop:int}).
Applying a twisted (co)homological technique 
we derive the linear Pfaffian system for $F_{L,N}$
(Theorem~\ref{thm:pfaff}).
In Sect.~\ref{sect:hgsol},
after a brief review of Lax formalism of ${\cal H}_{L,N}$ (Theorem~\ref{thm:lax}),
we solve it 
with special values of parameters.
Particular solutions thus obtained 
satisfy
 the same linear Pfaffian system as $F_{L,N}$.
This establishes a representation of the solutions in terms
of the hypergeometric functions (Theorem~\ref{thm:hgsol}).
 In the appendix we summarize the contiguity relations for $F_{L,N}$.

\section{Integral representation and Pfaffian system for hypergeometric function $F_{L,N}$}
\label{sect:hg}

In this section
we first introduce the integral representation of $F_{L,N}$,
then from which we derive 
the linear Pfaffian system 
by means of the viewpoint of twisted de Rham theory.
Conversely, the hypergeometric function $F_{L,N}$  
can be characterized as the unique
holomorphic solution of this Pfaffian system at
$0 \in {\mathbb C}^N$.

\subsection{Integral representations}

The hypergeometric function $F_{L,N}$, (\ref{eq:FLN}),
can be written as 
\begin{align}
\nonumber
&F_{L,N}(\alpha,\beta,\gamma;x)
= 
\prod_{k=1}^{L-1} \frac{\Gamma(\gamma_k)}{\Gamma(\alpha_k) \Gamma(\gamma_k-\alpha_k)}
\\
& \qquad
\times\int_{[0,1]^{L-1}}
\prod_{k=1}^{L-1} {z_k}^{\alpha_k-1}(1-z_k)^{\gamma_k-\alpha_k-1}
\prod_{i=1}^{N} (1-x_i z_1 z_2\cdots z_{L-1})^{-\beta_i}
{\rm d}z_1 \cdots {\rm d}z_{L-1},
\label{eq:rep_1}
\end{align}
provided $|x_i| <1$ and 
$ {\rm Re}(\gamma_k) > {\rm Re}(\alpha_k) >0$
and the branch of the integrand is 
assigned as
\[\arg z_k=\arg (1-z_k)=0 \quad \text{and}
\quad
\left|\arg (1-x_i z_1 z_2\cdots z_{L-1})\right| <\frac{\pi}{2}.
\]
Representation 
(\ref{eq:rep_1})
can be verified 
in a standard manner, i.e., 
by means of 
 the binomial theorem: 
\[ (1-z)^{-a}=\sum_{n=0}^\infty \frac{(a)_n}{(1)_n}z^n
\]
and the relation between the beta function and gamma function:
\[
B(a,b)=
\int_{0}^1 z^{a-1} (1-z)^{b-1} {\rm d}z
=\frac{\Gamma(a) \Gamma(b)}{\Gamma(a+b)},
\quad
{\rm Re}(a), {\rm Re}(b) >0.
\]
To transform the nonlinear form in the integrand of 
(\ref{eq:rep_1}) into a linear one,
we apply the change of integration variables 
\[
t_1=z_1, \quad t_2=z_1z_2, \quad
t_3=z_1z_2z_3, 
\quad  \ldots, \quad 
t_{L-1}=z_1z_2 \cdots z_{L-1}.
\]
The Jacobian of this transformation is 
calculated as 
\[
\left|
\frac{\partial(z_1,\ldots,z_{L-1})}{\partial(t_1,\ldots,t_{L-1})}
\right|
= \frac{1}{t_1t_2 \cdots t_{L-2}}.
\]
Put $t_0=1$ for convenience.
The integral representation of 
$F_{L,N}$
in which every factor of the integrand takes a linear form
is, therefore, established.

\begin{prop}  \label{prop:int}
Assume
$ {\rm Re} (\gamma_k) > {\rm Re} (\alpha_k) >0$.
For  $|x_i| <1$
it holds that
\begin{align}
\nonumber
&F_{L,N}(\alpha,\beta,\gamma;x)= 
\prod_{k=1}^{L-1} \frac{\Gamma(\gamma_k)}{\Gamma(\alpha_k) \Gamma(\gamma_k-\alpha_k)}
\\
& \qquad
\times\int_{\Delta}
{t_{L-1}}^{\alpha_{L-1}-1}
\prod_{k=1}^{L-2} {t_k}^{\alpha_k-\gamma_{k+1}}
\prod_{k=1}^{L-1} (t_{k-1}-t_k)^{\gamma_k-\alpha_k-1}
\prod_{i=1}^{N} (1-x_i t_{L-1})^{-\beta_i}
{\rm d}t_1 \cdots {\rm d}t_{L-1}
\label{eq:eular_rep}
\end{align}
with the integration domain $\Delta$ being an $(L-1)$-simplex
\[
\Delta  =
\{0 \leq t_{L-1} \leq \cdots \leq t_2 \leq t_1 \leq 1 \} \subset {\mathbb R}^{L-1}.
\]
Here 
the branch of the integrand is 
assigned as
\[\arg t_k=\arg (t_{k-1}-t_k)=0 \quad \text{and}
\quad
\left|\arg (1-x_i t_{L-1})\right| <\frac{\pi}{2}.
\]
\end{prop}

Based on this integral representation
and twisted de Rham theory,
we will discuss below the 
linear Pfaffian system characterizing 
$F_{L,N}$.

\subsection{Pfaffian system}

Consider a multivalued function 
\[
U(t)=
{t_{L-1}}^{\alpha_{L-1}}
\prod_{k=1}^{L-2} {t_k}^{\alpha_k-\gamma_{k+1}} 
\prod_{k=1}^{L-1} (t_{k-1}-t_k)^{\gamma_k-\alpha_k}
\prod_{i=1}^{N} (1-x_i t_{L-1})^{-\beta_i}
\]
with $t_0=1$
defined on 
\[
T
=\{t=(t_1,\ldots,t_{L-1}) \in {\mathbb C}^{L-1} \ | \  
t_k \neq 0, \ t_k \neq  t_{k-1}, \
t_{L-1} \neq 1/{x_i} \},
\]
which is the complement of singular locus 
$D= \bigcup_{k=1}^{L-1} \left(\{ t_k=0 \} \cup 
\{t_{k-1}-t_{k}=0\} \right)
\cup 
\bigcup_{i=1}^N\{ 1-x_i t_{L-1}=0\}$
of $U(t)$ in ${\mathbb C}^{L-1}$.
Let ${\cal L}$ be
the local system 
of rank one
determined by $1/U(t)$,
i.e.,
a flat line bundle consisting of 
the local solutions of 
$\nabla_\omega h =0$
on $T$,
where 
$\nabla_\omega$ is the covariant differential operator
given by
\begin{equation} \label{eq:cov}
\nabla_\omega={\rm d}+\omega \wedge,
\quad
\omega={\rm d} \log U(t).
\end{equation}
Let ${\cal L}^\vee$ 
be the dual local system of ${\cal L}$.
Denote by $H^p(T,{\cal L})$ (resp. $H_p(T,{\cal L}^\vee)$)
the $p$-th 
cohomology (resp. homology)
group with coefficients in ${\cal L}$ (resp. ${\cal L}^\vee$).
Under a certain genericity condition for the exponents
$\alpha_i,\beta_i,\gamma_i \in {\mathbb C} \setminus {\mathbb Z}$,
it holds that 
\[
\dim  H^p(T,{\cal L}) = \dim H_p(T,{\cal L}^\vee )
=  
\begin{cases}
N(L-1)+1 & \text{if $p=L-1$}
\\
0 & \text{if $p \neq L-1$}
\end{cases}
\]
and, furthermore,
bases of the top cohomology and homology groups 
are described as follows.
First we notice the isomorphism
\begin{align*}
H^{p}(T,{\cal L})
&\cong H^{p}(\Omega^\bullet(*D), \nabla_\omega)
\\
&=
\{
\xi \in \Omega^{p}(*D) \ | \ \nabla_\omega \xi =0
\}/ \nabla_\omega  \Omega^{p-1}(*D),
\end{align*}
where the right-hand side is
the de Rham cohomology group
determined by $\nabla_\omega$,
and 
$\Omega^p(*D)$
stands for the space of rational $p$-forms holomorphic outside $D$.
It can be verified that
the rational $(L-1)$-forms
\begin{align*}
\varphi_0&= \frac{{\rm d}t_1 \wedge \cdots \wedge {\rm d} t_{L-1}}{ t_{L-1} \prod_{k=1}^{L-1}(t_{k-1}-t_k)}, 
\\
\varphi_n^{(i)}&= \frac{{\rm d}t_1 \wedge \cdots \wedge {\rm d} t_{L-1}}{t_{L-1} (1-x_i t_{L-1})\prod_{\begin{subarray}{l} k=1 \\ k \neq n \end{subarray}}^{L-1}(t_{k-1}-t_k)}
\quad 
\left(
\begin{array}{c}
 1 \leq i \leq N \\
1 \leq n \leq L-1
\end{array}
\right)
\end{align*}
are cocycles
representing a basis of $H^{L-1}(T,{\cal L})$.
On the other hand,
a basis of 
$H_{L-1}(T,{\cal L}^\vee)$
can be constructed from the set of bounded chambers in the real locus 
$T \cap {\mathbb R}^{L-1}$
of $T$.
For simplicity, we fix the configuration $x \in {\mathbb C}^N$ of $N$ points 
to be real numbers such that
$0< x_N< \cdots <x_2<x_1<1$.
Accordingly, the set of bounded chambers 
is given by
\begin{align*}
{\Delta_0}' &= \{ 0< t_{L-1}< \cdots < t_1 < t_0 =1  \},
\\
{\Delta_n^{(i)}}'
&=
\left\{
\begin{array}{c} t_{L-1}>t_{L-2}>\cdots >t_{L-n-1}>0
\\
0<t_{L-n-1}< \cdots  < t_1 < t_0=1 \\ 
 1/x_{i-1}<t_{L-1}<1/x_i
\end{array}
\right\}
\quad 
\left(
\begin{array}{c}
 1 \leq i \leq N
 \\
1 \leq n \leq L-1
\end{array}
\right).
\end{align*}
The 
{\it regularizations} \cite{ak94, yos97}
of these cycles, denoted by
$\Delta_0$ and $\Delta_n^{(i)}$,
represent a basis of $H_{L-1}(T,{\cal L}^\vee)$.

Now
 we introduce the 
integrals
\begin{equation} \label{eq:integ}
y_0=\int_{\Delta} U(t) \varphi_0
\quad  \text{and}
\quad
y_n^{(i)}=\int_{\Delta} U(t) \varphi_n^{(i)}
\end{equation}
for any twisted cycle
$\Delta \in H_{L-1}(T,{\cal L}^\vee)$.
Then we have the

\begin{thm} \label{thm:pfaff}
The functions 
$y_0$
and $y_n^{(i)}$ satisfy the linear 
Pfaffian system 
\begin{subequations}
\label{subeq:pfaff}
\begin{align}
(x_i-1)
\frac{\partial y_0}{\partial x_i}
&= \beta_i 
\left(-y_0+
\sum_{m=1}^{L-1} y_m^{(i)} 
\right),
\\
(x_i-x_j) 
\frac{\partial y_n^{(j)}}{\partial x_i}
&= \beta_i
\left( y_n^{(i)}- y_n^{(j)}
\right), 
\\
x_i 
\frac{\partial y_n^{(i)}}{\partial x_i}
&=-\alpha_n y_n^{(i)}+(\gamma_n-\alpha_n) \sum_{m=n+1}^{L-1} y_m^{(i)}
+ \frac{\gamma_n-\alpha_n}{x_i-1}
\left(-y_0+
\sum_{m=1}^{L-1} y_m^{(i)} 
\right)
\nonumber
\\
&\quad
+ \sum_{\begin{subarray}{c} j=1 \\ j \neq i \end{subarray}}^N
\frac{\beta_j x_j}{x_i-x_j} 
\left( y_n^{(j)}- y_n^{(i)} \right).
\end{align}
\end{subequations}
\end{thm}

A proof of this theorem will be given in Sect.~\ref{subsec:proof}.

Let us consider the vector-valued function 
\[\vec{y}=\vec{y}(x; \Delta)=
{}^{\rm T}\left(
y_0,y_1^{(1)},y_2^{(1)},\ldots,y_{L-1}^{(1)},
y_1^{(2)},y_2^{(2)},\ldots,y_{L-1}^{(2)},
 \ldots,
y_1^{(N)},y_2^{(N)},\ldots,y_{L-1}^{(N)}
\right)
\]
defined by the integrals (\ref{eq:integ})
for $\Delta \in H_{L-1}(T,{\cal L}^\vee)$.
 Then (\ref{subeq:pfaff})
 takes the following expression:
\begin{equation} \label{eq:pfaff_matrix}
{\rm d} \vec{y}
= \left\{
\sum_{i=1}^N \left(E_i {\rm d} \log x_i+F_i  {\rm d} \log (x_i-1)
\right)+
\sum_{1 \leq i< j \leq N} G_{ij}   {\rm d} \log (x_i-x_j)
\right\}
\vec{y},
\end{equation}
where $E_i$, $F_i$ and $G_{ij}$ are
the square matrices of size $N(L-1)+1$:
\begin{allowdisplaybreaks}
\begin{align*}
&E_i=
\kbordermatrix{
&0 \hspace{-2mm}&&1&& 2 &&    &&&&i&&&&       &&\hspace{-2mm} N
\\
&& \vrule && \vrule && \vrule && \vrule &&&&&&\vrule&& \vrule
\\
\cline{2-18}
&  -a_1 \hspace{-2mm}&\vrule&  &\vrule&  &\vrule&  &\vrule& \hspace{-2mm}  b_{i,1}&&&&&\vrule&&\vrule& 
\\
& -a_2 \hspace{-2mm} &\vrule&  &\vrule&  &\vrule&  &\vrule& \hspace{-2mm}  a_2&b_{i,2} &&&&\vrule&&\vrule& 
\\
 i & -a_3 \hspace{-2mm} &\vrule& \hspace{-2mm} -\beta_1 I \hspace{-2mm} &\vrule& \hspace{-2mm} -\beta_2 I \hspace{-2mm}&\vrule& \hspace{-2mm}\cdots \hspace{-2mm} &\vrule& \hspace{-2mm} a_3&a_3& b_{i,3} &&&\vrule&\hspace{-2mm} \cdots \hspace{-2mm} &\vrule&\hspace{-2mm}  -\beta_N I 
 \\
 &\vdots \hspace{-2mm} &\vrule&  &\vrule&  &\vrule&  &\vrule&\hspace{-2mm}  \vdots&\vdots&\ddots&\ddots&&\vrule&&\vrule&  
 \\
 & -a_{L-1} \hspace{-2mm} &\vrule&  &\vrule&  &\vrule&  &\vrule& \hspace{-2mm} a_{L-1}& a_{L-1} &\cdots&a_{L-1}&b_{i,L-1} \hspace{-2mm}  &\vrule&&\vrule&  
\\
\cline{2-18}
&& \vrule && \vrule && \vrule && \vrule &&&&&&\vrule&& \vrule
},
\\
&F_i=
\kbordermatrix{
&0\hspace{-2mm}&& && & i & && \\
0&-\beta_i \hspace{-2mm} &\vrule& &\vrule& \hspace{-2mm}  \beta_i & \beta_i & \cdots& \beta_i \hspace{-2mm}  &\vrule& \\
\cline{2-11}
&&\vrule&&\vrule&&&&&\vrule& \\
\cline{2-11}
&a_1 \hspace{-2mm} &\vrule& &\vrule& \hspace{-2mm} -a_1 &-a_1& \cdots& -a_1 \hspace{-2mm} &\vrule& \\
i &a_2 \hspace{-2mm} &\vrule& &\vrule&\hspace{-2mm}  -a_2& -a_2&\cdots& -a_2 \hspace{-2mm}  &\vrule& \\
&\vdots\hspace{-2mm} &\vrule&&\vrule&\hspace{-2mm}  \vdots&\vdots&\ddots&\vdots \hspace{-2mm} &\vrule& \\
& a_{L-1} \hspace{-2mm} &\vrule& &\vrule&\hspace{-2mm}  -a_{L-1}& -a_{L-1}& \cdots&-a_{L-1}  \hspace{-2mm}  &\vrule& \\
\cline{2-11}
&&\vrule&&\vrule&&&&&\vrule& \\
},
\\
&G_{ij}=
\kbordermatrix{
&&& i &&&&j&& \\
&&\vrule&&\vrule&&\vrule&&\vrule&  \\
\cline{2-10}
i&& \vrule &\hspace{-2mm}  -\beta_j I \hspace{-2mm}  & \vrule&&\vrule&  \hspace{-2mm}  \beta_j I \hspace{-2mm}  &\vrule& \\
\cline{2-10}
&&\vrule&&\vrule&&\vrule&&\vrule& \\
\cline{2-10}
j&& \vrule & \hspace{-2mm}  \beta_i I  \hspace{-2mm}  & \vrule&&\vrule& \hspace{-2mm}  -\beta_i I \hspace{-2mm}  &\vrule& \\
\cline{2-10}
&&\vrule&&\vrule&&\vrule&&\vrule& \\
}, 
\end{align*}
\end{allowdisplaybreaks}
and 
 $a_n=\alpha_n-\gamma_n$ and 
$b_{i,n}=\sum_{j \neq i} \beta_j -\gamma_n $;  
the symbol $I$ denotes the identity matrix of size $L-1$.
We 
wrote a square matrix $M$ 
of size $N(L-1)+1$
as
\[
M=\kbordermatrix{
&0 \hspace{-2mm}&&1&&     && \hspace{-2mm}N \\
0&M_{00} \hspace{-2mm} & \vrule & \hspace{-2mm}  M_{01} \hspace{-2mm} &\vrule &\hspace{-2mm} \cdots \hspace{-2mm} &\vrule&
\hspace{-2mm}  M_{0N}
\\
\cline{2-8}
1& M_{10} \hspace{-2mm} & \vrule & \hspace{-2mm}  M_{11} \hspace{-2mm} &\vrule &\hspace{-2mm} \cdots \hspace{-2mm} &\vrule&
\hspace{-2mm} 
M_{1N}
\\
\cline{2-8}
      & \vdots \hspace{-2mm} &  \vrule & \hspace{-2mm}  \vdots\hspace{-2mm}  &\vrule &\hspace{-2mm}  \ddots \hspace{-2mm}  &\vrule&
\hspace{-2mm}  \vdots
\\
\cline{2-8}
N& M_{N0} \hspace{-2mm} & \vrule &\hspace{-2mm}  M_{N1}  \hspace{-2mm} &\vrule  &\hspace{-2mm}  \cdots \hspace{-2mm}  &\vrule&
\hspace{-2mm}  M_{NN}
}
\]
with dividing it into $(N+1)^2$ blocks 
so that $M_{00}$ becomes a scalar,
$M_{0j}$ $(j \neq 0)$ and $M_{i0}$ $(i \neq 0)$
row and  column $(L-1)$-vectors,
respectively,
and $M_{ij}$ $(i,j \neq 0)$ a square matrix of size $L-1$.

The linear Pffafian system (\ref{subeq:pfaff}), or (\ref{eq:pfaff_matrix}), is of rank $N(L-1)+1$
and the integrals 
\[
\vec{y}(x;\Delta_0) \quad \text{and} \quad \vec{y}(x;\Delta_n^{(i)})
\quad 
\left(
\begin{array}{c}
 1 \leq i \leq N
 \\
1 \leq n \leq L-1
\end{array}
\right)
\] 
provide a fundamental system of solutions.
In particular,
 $\vec{y}(x;\Delta_0)$ is the unique holomorphic solution
at
$0 \in {\mathbb C}^N$
up to multiplication by constants;
it is expressible in terms of the hypergeometric function
$F_{L,N}(\alpha,\beta,\gamma;x)$
according to the integral representation (see Proposition~\ref{prop:int})
as 
\begin{align*}
&y_0
= c F_{L,N}, \quad
y_1^{(i)}
= \frac{\gamma_1-\alpha_1}{\gamma_1}  
c F_{L,N}(\beta_i+1,\gamma_1+1),
\\
&y_2^{(i)}
= \frac{\alpha_1(\gamma_2-\alpha_2)}{\gamma_1\gamma_2} 
c F_{L,N}(\alpha_1+1,\beta_i+1,\gamma_1+1,\gamma_2+1),
\quad \ldots
\\
&y_n^{(i)}
=\frac{\alpha_1  \cdots \alpha_{n-1}(\gamma_n-\alpha_n)}{\gamma_1 \cdots \gamma_n} c F_{L,N}(\alpha_1+1,\ldots,\alpha_{n-1}+1,\beta_i+1,\gamma_1+1,\ldots,\gamma_{n}+1), 
\quad \ldots
\end{align*}
where $c=\prod_{k=1}^{L-1}\Gamma(\alpha_k)\Gamma(\gamma_k-\alpha_k)/\Gamma(\gamma_k)$. 
For notational simplicity,
we used the abbreviation 
$F_{L,N}(\beta_i+1, \gamma_1+1)$ 
to mean that among the parameters $(\alpha,\beta,\gamma)$
only the indicated ones $\beta_i$ and $\gamma_1$ are shifted by one,
and so forth.
We mention that
the differential equations satisfied by the first element
$y=y_0$ of $\vec{y}$ 
are indeed (\ref{eq:hgeq}).

\subsection{Verification of Theorem~\ref{thm:pfaff}}
\label{subsec:proof}

In general, it holds for an $(L-1)$-form $\varphi$  that
\[
\frac{\partial}{\partial x_i} \int_{\Delta} U \varphi
=\int_{\Delta}U
\left(\frac{1}{U}\frac{\partial U}{\partial x_i} \varphi + 
\frac{\partial \varphi}{\partial x_i}\right).
\]
Hence Theorem~\ref{thm:pfaff} 
is an immediate consequence of 
(\ref{eq:integ}) and the following lemma.

\begin{lemma}
Define a linear operator 
$\nabla_i$ $(i=1,\ldots,N)$
acting on a differential form 
$\varphi$
by
\[
\nabla_i \varphi= \frac{1}{U}\frac{\partial U}{\partial x_i} \varphi + 
\frac{\partial \varphi}{\partial x_i}.
\]
The rational $(L-1)$-forms
$\varphi_0$ and $\varphi_n^{(i)}$
satisfy the relations
\begin{subequations}
\begin{align} \label{eq:form1}
(x_i-1)\nabla_i \varphi_0
&= \beta_i 
\left(
 - \varphi_0
+ \sum_{m=1}^{L-1} \varphi_m^{(i)} 
\right),
\\
 \label{eq:form2}
(x_i-x_j) \nabla_i \varphi_n^{(j)}
&= \beta_i
\left( \varphi_n^{(i)}- \varphi_n^{(j)}
\right), 
\\
 \label{eq:form3}
x_i \nabla_i \varphi_n^{(i)}
&\equiv
-\alpha_n \varphi_n^{(i)}+(\gamma_n-\alpha_n) \sum_{m=n+1}^{L-1} \varphi_m^{(i)}
+ \frac{\gamma_n-\alpha_n}{x_i-1}
\left(
 - \varphi_0 +
\sum_{m=1}^{L-1} \varphi_m^{(i)} 
\right)
\nonumber
\\
&\quad
+ \sum_{\begin{subarray}{c} j=1 \\ j \neq i \end{subarray}}^N
\frac{\beta_j x_j}{x_i-x_j} 
\left( \varphi_n^{(j)}- \varphi_n^{(i)} \right)
\quad 
(\text{modulo $\nabla_\omega \Omega^{L-2}(*D)$}).
\end{align}
\end{subequations}
\end{lemma}

\pf
We will use the notation
\begin{align*}
\underline{{\rm d} t}&= {\rm d}t_1 \wedge \cdots \wedge {\rm d} t_{L-1},
\\
*{\rm d}t_j&=(-1)^{j-1} {\rm d} t_1 \wedge \cdots  \wedge
 \widehat{ {\rm d} t_{j} }  \wedge \cdots  \wedge {\rm d} t_{L-1};
\end{align*}
therefore, ${\rm d}t_j \wedge *{\rm d}t_j= \underline{{\rm d} t}$.
We abbreviate 
$\prod_{k=1}^{L-1}$ and $\prod^{L-1}_{\begin{subarray}{c}k=1 \\ k \neq n \end{subarray} }$ respectively
to $\prod_k$ and $\prod_{k \neq n}$,
and so forth.
From the definition 
\begin{equation} \label{eq:help0}
\varphi_0= \frac{\underline{{\rm d} t}}{ t_{L-1} \prod_{k}(t_{k-1}-t_k)}, 
\quad
\varphi_n^{(i)}= \frac{\underline{{\rm d} t}}{t_{L-1} (1-x_i t_{L-1})
\prod_{k \neq n} (t_{k-1}-t_k)},
\end{equation}
it is readily seen that
\begin{align} \label{eq:help__1}
\sum_{m=n+1}^{L-1} \varphi_m^{(i)}&=
\frac{ t_n-t_{L-1} }{ t_{L-1}(1-x_i t_{L-1}) \prod_k (t_{k-1}-t_k)}\underline{{\rm d} t}
\quad (0 \leq  n \leq L-2),
\\
\label{eq:help2}
-\varphi_0 +\sum_{m=1}^{L-1} \varphi_m^{(i)}
&=
\frac{ x_i-1}{ (1-x_i t_{L-1}) \prod_k (t_{k-1}-t_k)}
\underline{{\rm d} t}.
\end{align} 
Since $\varphi_0$ does not depend on $x_i$
it follows that
\[
\nabla_i \varphi_0 
=\frac{1}{U}\frac{\partial U}{\partial x_i}  \varphi_0
= \frac{\beta_i t_{L-1}}{1-x_i t_{L-1}}  \varphi_0
= \frac{  \beta_i  \underline{{\rm d}t}}{ (1-x_i t_{L-1}) \prod_k (t_{k-1}-t_k)},
\]
which coincides with (\ref{eq:form1}) 
according to (\ref{eq:help2}).
Likewise (\ref{eq:form2}) can be verified as 
\begin{align*}
\nabla_{i} \varphi_n^{(j)}
&= \frac{ \beta_i \underline{{\rm d}t} }{ (1-x_i t_{L-1}) (1-x_j t_{L-1}) 
\prod_{k \neq n}(t_{k-1}-t_k)}
\\
&=\frac{\beta_i}{x_i-x_j} \left(\frac{1}{1-x_i t_{L-1}}-\frac{1}{1-x_j t_{L-1}}\right)
\frac{\underline{{\rm d}t} }{ t_{L-1}
\prod_{k \neq n} (t_{k-1}-t_k)}
\\
&= \frac{\beta_i}{x_i-x_j} \left( \varphi_n^{(i)}-\varphi_n^{(j)} \right).
\end{align*}

It is difficult to calculate directly $\nabla_i \varphi_n^{(i)}$ 
because
$\varphi_n^{(i)}$ depends on $x_i$.
So we first prepare an appropriate coboundary,
by which we eliminate the $x_i$ dependence of
$\varphi_n^{(i)}$
modulo $\nabla_\omega \Omega^{L-2}(*D)$. 
Consider the rational $(L-2)$-form
\begin{align*}
\Omega^{L-2}(*D) \ni
\xi_n
&=
\prod_{k \neq n} {\rm d} \log (t_{k-1}-t_k) \wedge 
\\
&=\frac{(-1)^{L+n-1}}{ \prod_{k \neq n} (t_{k-1}-t_k) } \sum_{j=n}^{L-1} 
*{\rm d}t_j.
\end{align*}
Its covariant derivative reads
(recall (\ref{eq:cov}))
\begin{align*}
\nabla_\omega \xi_n
&=
(-1)^{L+n}
\left(
\frac{\gamma_n-\alpha_n}{ \prod_k (t_{k-1}-t_k)}
-\frac{\alpha_{L-1}}{t_{L-1} \prod_{k \neq n} (t_{k-1}-t_k)}
\right.
\\
& \qquad  \left.
-\sum_{j=1}^{N} 
\frac{ \beta_j x_j}{(1-x_j t_{L-1}) \prod_{k \neq n} (t_{k-1}-t_k) }
+
\sum_{j=n}^{L-2} \frac{\gamma_{j+1}-\alpha_{j}}{t_{j} \prod_{k \neq n}(t_{k-1}-t_k)}
\right)
\underline{ {\rm d} t}.
\end{align*}
Hence
\begin{align*}
\sum_{j=1}^N \beta_j \varphi_n^{(j)} 
&=
\sum_{j=1}^N \beta_j
\left(
\frac{1}{t_{L-1} \prod_{k \neq n} (t_{k-1} -t_k ) }
+ \frac{x_j}{(1-x_j t_{L-1}) \prod_{k \neq n} (t_{k-1} -t_k ) }
\right) \underline{ {\rm d} t}
\\
&\equiv
\left(
\frac{\gamma_n-\alpha_n}{ \prod_k (t_{k-1}-t_k)}
+
\frac{\sum_{j=1}^N \beta_j-\alpha_{L-1}}{t_{L-1} \prod_{k \neq n} (t_{k-1}-t_k)}
+
\sum_{j=n}^{L-2} \frac{\gamma_{j+1}-\alpha_{j}}{t_{j} \prod_{k \neq n}(t_{k-1}-t_k)}
\right)
 \underline{ {\rm d} t}
\end{align*}
modulo $\nabla_\omega \xi_n$. 
Now the derivative can be calculated as
\begin{align*}
\left(\frac{x_i}{\beta_i} \nabla_i+1\right)
\sum_{j=1}^N \beta_j \varphi_n^{(j)} 
&=\left( \frac{1}{1-x_i t_{L-1}} + \frac{x_i}{\beta_i} \frac{\partial}{\partial x_i}\right)\sum_{j=1}^N \beta_j \varphi_n^{(j)} 
\\
&\equiv
\left(
\frac{\gamma_n-\alpha_n}{ (1-x_i t_{L-1})\prod_k (t_{k-1}-t_k)}
+
\frac{\sum_{j=1}^N \beta_j-\alpha_{L-1}}{t_{L-1} (1-x_i t_{L-1})\prod_{k \neq n} (t_{k-1}-t_k)}
\right.
 \\
& \qquad \qquad
\left.
+
\sum_{j=n}^{L-2} \frac{\gamma_{j+1}-\alpha_{j}}{t_{j} (1-x_i t_{L-1}) \prod_{k \neq n}(t_{k-1}-t_k)}
\right)
 \underline{ {\rm d} t}.
\end{align*}
Applying
(\ref{eq:help0}), (\ref{eq:help2}) and 
Lemma~\ref{lemma:sub} below,
we thus arrive at
\[
\left(\frac{x_i}{\beta_i} \nabla_i+1\right)
\sum_{j=1}^N \beta_j \varphi_n^{(j)} 
\equiv
\frac{\gamma_n-\alpha_n}{x_i-1} 
 \left( -\varphi_0+\sum_{m=1}^{L-1} \varphi_m^{(i)}  \right)
 + \left( \sum_{j=1}^N \beta_j - \alpha_{n} \right) 
 \varphi_n^{(i)} 
 +(\gamma_n-\alpha_n) \sum_{m=n+1}^{L-1} \varphi_m^{(i)},
\]
which 
establishes 
(\ref{eq:form3}) by virtue of (\ref{eq:form2}).
The proof of the lemma is complete.
\qed

\begin{lemma} \label{lemma:sub}
One has
\begin{equation}  \label{eq:claim}
\sum_{j=n}^{L-2} \frac{\gamma_{j+1}-\alpha_{j}}{t_{j} (1-x_i t_{L-1}) \prod_{k \neq n}(t_{k-1}-t_k)}
 \underline{ {\rm d} t} \equiv 
 (\alpha_{L-1}-\alpha_n) \varphi_n^{(i)}
 +(\gamma_n-\alpha_n) \sum_{m=n+1}^{L-1} \varphi_m^{(i)}
 \end{equation}
 modulo $\nabla_\omega \Omega^{L-2}(*D)$. 
\end{lemma}

\pf 
Taking partial fraction decomposition 
yields
\begin{align}
\lefteqn{\sum_{j=n}^{L-2} \frac{\gamma_{j+1}-\alpha_{j}}{t_{j} (1-x_i t_{L-1}) \prod_{k \neq n}(t_{k-1}-t_k)}
 \underline{ {\rm d} t}
  }
\nonumber  \\
 &
 =\frac{1}{t_{L-1}(1-x_i t_{L-1})}
 \sum_{j=n}^{L-2} (\gamma_{j+1}-\alpha_j)
 \left(
 \frac{1}{\prod_{k \neq n} (t_{k-1}-t_k)}
 - \frac{1}{t_j} \sum_{m=j+1}^{L-1}   \frac{1}{\prod_{k \neq n,m} (t_{k-1}-t_k)}
 \right)
 \underline{ {\rm d} t}
 \nonumber
 \\
 &=  
 \sum_{j=n}^{L-2} (\gamma_{j+1}-\alpha_j)  \varphi_n^{(i)}
 -
 \frac{1}{t_{L-1}(1-x_i t_{L-1})}
\sum_{j=n}^{L-2} \sum_{m=j+1}^{L-1}
 \frac{\gamma_{j+1}-\alpha_{j}}{t_{j} \prod_{k \neq n, m}(t_{k-1}-t_k)} \underline{ {\rm d} t}.
 \label{eq:pfd}
 \end{align}
On the other hand, 
we consider 
the $(L-2)$-form
\begin{align*}
\Omega^{L-2}(*D) \ni \psi_{n,m}
&=\left(\prod_{k \neq n,m} {\rm d} \log (t_{k-1}-t_k) \wedge \right)
{\rm d} t_{L-1}
\\
&=\frac{(-1)^{n+m-1}}{ \prod_{k \neq n,m} (t_{k-1}-t_k) } \sum_{j=n}^{m-1} * {\rm d}t_j
\end{align*}
for $n < m$,
of which the covariant derivative reads
\[
\nabla_\omega \psi_{n,m}=
(-1)^{n+m}
\left(\sum_{j=n}^{m-1} \frac{\gamma_{j+1}-\alpha_{j}}{t_{j} \prod_{k \neq n, m}(t_{k-1}-t_k)}
+ \frac{\gamma_n-\alpha_n}{ \prod_{k \neq m}(t_{k-1}-t_k) }
- \frac{\gamma_m-\alpha_m}{ \prod_{k \neq n}(t_{k-1}-t_k) }
\right) 
\underline{ {\rm d} t}.
\]
Observe that
$\nabla_\omega \psi_{m,n}$ still remains a coboundary
if we multiply it by any rational function $g \in \Omega^0(*D)$
such that 
$\partial g/\partial t_k =0$ $(\forall k \neq L-1)$.
In particular, 
by choosing $g={t_{L-1}}^{-1}(1-x_i t_{L-1})^{-1}$,
we have
\[
\frac{1}{t_{L-1}(1-x_i t_{L-1})}
\sum_{j=n}^{m-1} \frac{\gamma_{j+1}-\alpha_{j}}{t_{j} \prod_{k \neq n, m}(t_{k-1}-t_k)} \underline{ {\rm d} t}
+(\gamma_n-\alpha_n)\varphi_m^{(i)}
-(\gamma_m-\alpha_m)\varphi_n^{(i)}
\equiv 0.
\]
Summation over $m=n+1, \ldots, L-1$ of this formula
entails
\[
\frac{1}{t_{L-1}(1-x_i t_{L-1})}
\sum_{m=n+1}^{L-1}
\sum_{j=n}^{m-1} \frac{\gamma_{j+1}-\alpha_{j}}{t_{j} \prod_{k \neq n, m}(t_{k-1}-t_k)} \underline{ {\rm d} t}
+(\gamma_n-\alpha_n) \sum_{m=n+1}^{L-1}\varphi_m^{(i)}
-\sum_{m=n+1}^{L-1}(\gamma_m-\alpha_m)\varphi_n^{(i)}
\equiv 0.
\]
Substituting the above into (\ref{eq:pfd}),
we verify the desired result (\ref{eq:claim}).
\qed

\section{Hypergeometric solution of Hamiltonian system ${\cal H}_{L,N}$}
\label{sect:hgsol}

In this section we first review Lax formalism of ${\cal H}_{L,N}$ 
following \cite{tsu10}.
Under a certain condition of parameters 
${\cal H}_{L,N}$ admits
particular solutions 
such that the associated Fuchsian system becomes reducible.
These solutions are
governed by the Pfaffian system derived in 
the previous section and, 
thereby, expressible in terms of the hypergeometric function.

\subsection{Lax formalism of ${\cal H}_{L.N}$}

We begin with a brief review of Lax formalism of ${\cal H}_{L,N}$.
Consider an $L \times L$ Fuchsian system 
\begin{align}  \label{eq:A}
\frac{\partial \Phi}{\partial z} &= A \Phi = \sum_{i=0}^{N+1} \frac{A_i}{z-u_i} \Phi
\end{align}
with $N+3$ regular singularities 
$\{u_0=1, u_1, \ldots, u_{N}, u_{N+1}=0, u_{N+2}=\infty\} \subset {\mathbb P}^1$,
of which the characteristic exponents
at each singularity $z=u_i$,
i.e., the eigenvalues of each residue matrix $A_i$,
are listed in the following table (Riemann scheme):
\[
\begin{array}{cc}
\hline
\text{Singularity} & \text{Exponents}  \\ \hline
u_i  \mbox{\ } (0 \leq i \leq N )& (-\theta_i, 0, \ldots,0)  \\
u_{N+1}=0 & (e_0,e_1,\ldots,e_{L-1})  \\     
u_{N+2}=\infty & (\kappa_{0}-e_0,\kappa_{1}-e_1, \ldots,\kappa_{L-1}-e_{L-1} )  \\ 
\hline
   \end{array}
\]
We can, and will, normalize the exponents as
${\rm tr \,} A_{N+1}=\sum_{n=0}^{L-1} e_n={(L-1)}/{2}$
without loss of generality.
Assume 
$\sum_{n=0}^{L-1}\kappa_n
=\sum_{i=0}^N \theta_i$
(Fuchs' relation).
Such Fuchsian systems as above 
then turn out to
constitute a $2N(L-1)$-dimensional family
and, actually, can be written
in terms of the accessory parameters
$b_n^{(i)}$ and $c_n^{(i)}$ 
in the following way:
\begin{align*}
A_i 
&=  {}^{\rm T}\left(b_0^{(i)}, b_1^{(i)}, \ldots ,b_{L-1}^{(i)}\right)
\cdot
\left(c_0^{(i)}, c_1^{(i)}, \ldots ,c_{L-1}^{(i)}\right)\quad
(0 \leq i \leq N),
\\
A_{N+1} 
&= 
 \begin{pmatrix}
   e_0 &   w_{0,1}      & \cdots  & w_{0,L-1}
   \\
        & e_1 &  \ddots       & \vdots
   \\     
          &        &\ddots& w_{L-2,L-1}
   \\
          &        &          & e_{L-1}    
\end{pmatrix}, 
\end{align*}
where
$c_0^{(i)}=1$ and
 $w_{m,n}=-\sum_{i=0}^N b_m^{(i)} c_n^{(i)}$.
We thus find the relations
\begin{equation} \label{eq:k&t}
\left({\rm tr \,} A_i
 =\right)
  \sum_{n=0}^{L-1} b_n^{(i)}c_n^{(i)}
 = -\theta_i
 \quad  \text{and} \quad
 \sum_{i=0}^{N} b_n^{(i)} c_n^{(i)}=- \kappa_n,
\end{equation}
the latter of which comes from the diagonal entries of 
 the lower triangular matrix  
 $A_{N+2}=-\sum_{i=0}^{N+1}A_i$.
Since 
$A_{N+1}$ and $A_{N+2}$ 
are triangular,  
there still remains the degree of freedom of
a similarity transformation by a diagonal matrix.
Consequently, 
the essential number of the
accessory parameters 
is confirmed to be
$2N(L-1)$.
In fact, they can be realized 
by the canonical variables
$q_n^{(i)}$ and $p_n^{(i)}$
of ${\cal H}_{L,N}$; see (\ref{eq:qp}) below.

The isomonodromic family of Fuchsian systems
of the form (\ref{eq:A}) is described by 
the  integrability condition of the extended linear system, i.e.,
 (\ref{eq:A}) itself and its deformation equations 
\begin{equation}  \label{eq:B}
\frac{\partial \Phi}{\partial u_i} 
=  B_i \Phi, \quad
B_i
=\frac{A_i}{u_i-z} - \frac{1}{u_i} 
\begin{pmatrix} 
-\frac{\theta_i}{L} &
\\
& \ddots &
\\
\hsymb{*}&&-\frac{\theta_i}{L}
\end{pmatrix}
\quad (1 \leq i \leq N).
\end{equation}
Here the lower triangular part ($*$)
of the second term is exactly the same as $A_i$.

\begin{thm}[See \cite{tsu10}]  \label{thm:lax}
The integrability condition 
\begin{equation} \label{eq:lax}
 \frac{\partial A}{\partial u_i} -\frac{\partial B_i}{\partial z}+[A,B_i]=0
\end{equation}
 of {\rm(\ref{eq:A})} and {\rm(\ref{eq:B})}
is equivalent to the polynomial Hamiltonian system
${\cal H}_{L,N}$, {\rm(\ref{subeq:HLN})},
via the change of variables
\begin{equation}  \label{eq:qp}
x_i=\frac{1}{u_i}, \quad
q_n^{(i)}=\frac{c_n^{(i)}}{c_n^{(0)}} \quad \text{and} \quad 
p_n^{(i)}=-b_n^{(i)}c_n^{(0)}.
\end{equation}
\end{thm}

\subsection{Particular solution of ${\cal H}_{L,N}$}
Suppose that $\kappa_0-\sum_{i=1}^N \theta_i=0$.
This condition enables us to
restrict
$b_n^{(i)}$ and $c_n^{(i)}$
to the subvariety
\begin{align*}
&b_0^{(0)}=0, \quad b_0^{(i)}=-\theta_i \quad (i \neq 0),
\\
&
c_1^{(0)}=\cdots =c_{L-1}^{(0)}
\quad
\text{and} \quad
c_n^{(i)}=0 \quad 
(i \neq 0,n \neq 0)
\end{align*}
while keeping consistency of the linear system
(\ref{eq:A}) and (\ref{eq:B}).
Notice here that this restriction amounts to $q_n^{(i)}=0$.
In view of
(\ref{eq:k&t}) we see that
the matrices $A_i$ and $B_i$ can be 
parameterized 
by the $N(L-1)+1$ variables
$f:=
1/c_1^{(0)}=\cdots =1/c_{L-1}^{(0)}$ and $b_n^{(i)}$ $(i \neq 0, n \neq 0)$.
It actually follows that 
\begin{align*}
A_0 &= 
 \begin{pmatrix}
 0 & 0 & \cdots & 0
 \\
- \kappa_1 f& - \kappa_1 & \cdots  & - \kappa_1
\\
\vdots & \vdots & \ddots & \vdots
\\
- \kappa_{L-1} f & - \kappa_{L-1} & \cdots  & - \kappa_{L-1}
\end{pmatrix}, 
\quad
A_i =
\begin{pmatrix}
 -\theta_i & 0 & \cdots & 0
 \\
b_1^{(i)}& 0 & \cdots  & 0 
\\
\vdots & \vdots & \ddots & \vdots
\\
b_{L-1}^{(i)} & 0& \cdots  & 0
\end{pmatrix}
\quad (1 \leq i \leq N),
\\
A_{N+1} 
&= 
 \begin{pmatrix}
   e_0 &   0   & 0  & 0 & \cdots  & 0
   \\
        & e_1 & \kappa_1 & \kappa_1 &  \cdots     &  \kappa_1
   \\     
          &        &e_2& \kappa_2 &  \cdots     &  \kappa_2
   \\
          &        &          & e_3 &  \ddots     &  \vdots
   \\
          &         &          &           & \ddots  & \kappa_{L-2}
    \\
          &         &          &            &              & e_{L-1}
\end{pmatrix}, 
\end{align*}
therefore, (\ref{eq:A}) is clearly reducible.
On the other hand, we have 
\[
B_i=
\frac{\theta_i}{L u_i} 
\begin{pmatrix}
 1-L & & & 
 \\
&1&&
\\
&&\ddots&
\\
&&& 1
\end{pmatrix}
+
\frac{z}{u_i(u_i-z)}
\begin{pmatrix}
 -\theta_i & 0 & \cdots & 0
 \\
b_1^{(i)}& 0 & \cdots  & 0 
\\
\vdots & \vdots & \ddots & \vdots
\\
b_{L-1}^{(i)} & 0& \cdots  & 0
\end{pmatrix}.
\]
Observe that in this situation
only the $(n,0)$-entries $(1 \leq n \leq L-1)$
of the matrix equation (\ref{eq:lax}) are nontrivial.
In fact, 
\begin{align*}
\left( \frac{\partial A}{\partial u_i} -\frac{\partial B_i}{\partial z} \right)_{n,0}
&=
\frac{-\kappa_n}{z-1} \frac{\partial f}{\partial u_i}
+ \sum_{j=1}^N \frac{1}{z-u_j} \frac{\partial b_n^{(j)}}{\partial u_i},
\\
[A,B_i]_{n,0}
&=
\frac{\kappa_n }{u_i(z-1)(z-u_i)}  
\left( -\theta_i u_i f +
z \sum_{m=1}^{L-1} b_m^{(i)}
\right)
\\
& \quad
+ \frac{1}{u_i(z-u_i)} 
\left(  
(e_0-e_n)b_n^{(i)}  - \kappa_n \sum_{m=n+1}^{L-1} b_m^{(i)}
\right) 
\\
& \quad
+ \frac{1}{u_i(z-u_i)}\sum_{j=1}^N
\frac{\theta_i b_n^{(j)} u_i- \theta_j b_n^{(i)}z}{z-u_j}.
\end{align*}
Residue calculus at $z=1$, 
$z=u_j$ $(j \neq i)$ and $z=u_i$ 
yields the system of
 linear differential equations
for unknowns $f$ and $b_n^{(i)}$:
\begin{subequations} \label{subeq:eq_f}
\begin{align} 
u_i(1-u_i) \frac{\partial f}{\partial u_i}
&=-u_i \theta_i f + \sum_{m=1}^{L-1}  b_m^{(i)},
\\
(u_i-u_j) \frac{\partial b_n^{(j)}}{\partial u_i}
&= \theta_i b_n^{(j)} -\frac{u_j}{u_i} \theta_j b_n^{(i)},
\\
\nonumber
\frac{\partial b_n^{(i)}}{\partial u_i}
&=
\frac{1}{u_i} \left(
(e_n-e_0+\theta_i) b_n^{(i)}+\kappa_n \sum_{m=n+1}^{L-1} b_m^{(i)}
\right)
+ \frac{\kappa_n}{u_i-1} 
\left(  \theta_i f- \sum_{m=1}^{L-1} b_m^{(i)} \right)
\\
&\quad
- \sum_{\begin{subarray}{l} j=1 \\ j \neq i \end{subarray}}^N
\frac{\theta_i b_n^{(j)} - \theta_j b_n^{(i)}}{u_i-u_j}. 
\end{align} 
\end{subequations}
If we apply the change of variables
\[
x_i=\frac{1}{u_i},  \quad y_0 = \frac{f}{\prod_{j=1}^N {u_j}^{\theta_j}}
\quad
\text{and} \quad
y_n^{(i)}= \frac{b_n^{(i)}}{\theta_i \prod_{j=1}^N {u_j}^{\theta_j} },
\]
then (\ref{subeq:eq_f})
is converted into the Pfaffian system
 for $F_{L,N}$
 (see (\ref{subeq:pfaff}) in Theorem~\ref{thm:pfaff})
with $\alpha_n=e_n-e_0$,  $\beta_n=  - \theta_n$ and  
$\gamma_n= e_n-e_0-\kappa_n$.
Combining this fact with Theorem~\ref{thm:lax},
we finally arrive at the

\begin{thm} \label{thm:hgsol}
When $\kappa_0-\sum_{i=1}^N \theta_i=0$,
the Hamiltonian system ${\cal H}_{L,N}$ 
admits
a particular solution 
of the form
\[
q_n^{(i)} = 0, \quad
p_n^{(i)}=-\theta_i\frac{y_n^{(i)}}{y_0}
\quad 
\left(
\begin{array}{c}
 1 \leq i \leq N
 \\
1 \leq n \leq L-1
\end{array}
\right)
\]
where $\{y_0,y_n^{(i)}\}$ is an arbitrary solution of the linear Pfaffian system, {\rm(\ref{subeq:pfaff})} or {\rm(\ref{eq:pfaff_matrix})},
with
\[
\alpha_n=e_n-e_0, 
\quad 
\beta_n=  - \theta_n, \quad 
\gamma_n= e_n-e_0-\kappa_n.
\]
\end{thm}

\begin{remark}\rm
We already know that
$y_0$ is a solution of the hypergeometric equation
(\ref{eq:hgeq}).
Moreover, it is possible to write all the other elements $y_n^{(i)}$ 
as linear combinations of derivatives of $y_0$.
In fact, we can carry out it by the differential operators
appearing in  the contiguity relations for the 
hypergeometric series
$F_{L,N}$;
cf. Sect.~\ref{sect:hg} and the appendix below.
\end{remark}

\begin{remark}\rm
Particular solutions of the Schlesinger system
for the case of a general spectral type have been
studied by Dubrovin--Mazzocco \cite{dm07}.
For instance
they showed that
if the monodromy group of the associated Fuchsian system is triangular, 
then the Schlesinger system can be solved in terms of 
solutions of linear differential equations
(that are generally inhomogeneous).
\end{remark}

\appendix
\section{Contiguity relations for $F_{L,N}$}

In this appendix we provide a table of 
the contiguity relations for 
$F=F_{L,N}(\alpha,\beta,\gamma;x)$.
We shall use again such an abbreviation as 
$F(\alpha_n+1)$ to represent the same function as
$F$ except increasing the indicated parameter
$\alpha_n$ by one.
Recall the notation 
(\ref{eq:eular})
of the Eular operators.

\begin{thm}  \label{thm:cont}
The hypergeometric function $F=F_{L,N}(\alpha,\beta,\gamma;x)$
satisfies the contiguity relations
\begin{align}
&F(\alpha_n+1)=\frac{{\cal D}+\alpha_n}{\alpha_n} F,
\label{eq:cont1}
\\
&F(\beta_i+1)=\frac{\delta_i+\beta_i}{\beta_i} F, 
\label{eq:cont2}
\\
&F(\gamma_n+1)=
\frac{\gamma_n}{\varepsilon_L}
\left\{
\sum_{i=1}^N \frac{\partial }{\partial x_i}
\prod^{L-1}_{\begin{subarray}{l} k=1 \\ k \neq n\end{subarray}} 
({\cal D}+\gamma_k-1) 
- \sum_{j=0}^{L-1} \varepsilon_j ({\cal D}+\gamma_n)^{L-1-j}
\right\}F,
\label{eq:cont3}
\\
&F(\alpha_n-1)=
\frac{\alpha_n-1}{\varepsilon'_L}
\left\{
\sum_{i=1}^N x_i(\delta_i+\beta_i)
\prod^{L-1}_{\begin{subarray}{l} k=1 \\ k \neq n\end{subarray}} 
({\cal D}+\alpha_k) 
- \sum_{j=0}^{L-1} \varepsilon'_j ({\cal D}+\alpha_n-1)^{L-1-j}
\right\}F,
\label{eq:cont4}
\\
&F(\gamma_n-1)=\frac{{\cal D}+\gamma_n-1}{\gamma_n-1} F,
\label{eq:cont5}
\\
&F(\beta_i+1,\beta_j-1)=
\frac{(x_i-x_j)\frac{\partial}{\partial x_i} + \beta_i}{\beta_i} F,
\label{eq:cont6}
\\
&F(\alpha_1+1,\ldots,\alpha_{L-1}+1,\beta_i+1,
\gamma_1+1,\ldots,\gamma_{L-1}+1)
=  \frac{\gamma_1\cdots\gamma_{L-1}}{\alpha_1 \cdots \alpha_{L-1} \beta_i} \frac{\partial F}{\partial x_i}.
\label{eq:cont7}
\end{align}
Here each
$\varepsilon_j$ 
denotes the $j$-th elementary symmetric polynomial in 
$L$ variables
$\alpha_k-\gamma_n$
$(k=1,\ldots,L-1)$
and
$\sum_{i=1}^N \beta_i -\gamma_n$,
and similarly
$\varepsilon'_j$ 
does 
that in 
$\gamma_k-\alpha_n$
$(k=1,\ldots,L-1)$
and
$1-\alpha_n$.
\end{thm}

\pf
Notice the formula
$\delta_i x_i=x_i(\delta_i+1)$
with $x_i$ regarded as the operator multiplying $x_i$.
Then (\ref{eq:cont1}), (\ref{eq:cont2}), (\ref{eq:cont5})
and (\ref{eq:cont7})
are immediate
from the definition (\ref{eq:FLN}) of the hypergeometric series.
By an analogue of
the classical factorization method, 
(\ref{eq:cont3}), (\ref{eq:cont4}) and  (\ref{eq:cont6}) 
can be obtained.
\qed
\\

Let ${\cal S}={\cal S}(\alpha,\beta,\gamma)$
be the linear space of solutions of 
the hypergeometric equation (\ref{eq:hgeq}).  
In general, the linear operators
appearing in Theorem~\ref{thm:cont}
induce isomorphisms of these spaces.
For instance, let $H$ and $B$
be the differential operators
defined by
\[
H =
\sum_{i=1}^N \frac{\partial }{\partial x_i}
\prod^{L-1}_{\begin{subarray}{l} k=1 \\ k \neq n\end{subarray}} 
({\cal D}+\gamma_k-1) 
- \sum_{j=0}^{L-1} \varepsilon_j ({\cal D}+\gamma_n)^{L-1-j},
\quad
B ={\cal D}+\gamma_n;
\]
cf.  (\ref{eq:cont3}) and (\ref{eq:cont5}).
The linear homomorphisms
\[
H : {\cal S} \to {\cal S}(\gamma_n+1),
\quad
B  :{\cal S}(\gamma_n+1) \to {\cal S}
\]
are isomorphisms if and only if 
$\gamma_n \varepsilon_{L} \neq 0$.

\small
\paragraph{\it Acknowledgement.}
The author 
is grateful to
Masaaki Yoshida
for reading carefully 
the first draft of this article
and giving helpful comments.
He has also benefited from discussions with
Yousuke Ohyama and Yasuhiko Yamada.


\end{document}